\documentclass[12pt, final]{amsart}
\usepackage[cp850]{inputenc}
\usepackage{amssymb}
\usepackage[mathscr]{eucal}
\usepackage{amscd}
\newcommand{\R}{\mathbb R}

\def\Ext{\operatorname{Ext}}

\theoremstyle{plain}

\newtheorem{thm}{Theorem}

\newtheorem{proposition}{Proposition}
\newtheorem{lemma}{Lemma}

\theoremstyle{remark}

\title{Extensions by spaces of continuous functions}
\thanks{Supported in part by DGICYT project 2001-0813.}
\thanks{2000 Math. Subject Class.  46B03, 46B07}
\author{Jes\'us M. F. Castillo}
\author{Yolanda Moreno}
\address{Departamento de Matem\'aticas, Universidad de Extremadura, Avenida
de Elvas, 06071 Badajoz, Spain}

\email{castillo@unex.es} \email{ymoreno@unex.es}
\begin{document}
\maketitle
\begin{abstract} We characterize the Banach spaces $X$ such that $\Ext(X, C(K))=0$ for every compact space
$K$.
\end{abstract}\medskip

\section{Introduction}

In this paper we give a necessary and sufficient condition on a
Banach space $X$ so that every The exact sequence$$\begin{CD}
0@>>>C(K) @>>> E@>>> X @>>> 0,
\end{CD}$$ splits. Such exact sequences are sometimes called
\emph{extensions of $X$ by $C(K)$}. Thus, we  shall characterize
those Banach spaces $X$ such that, for every compact Hausdorff
space $K$ one has $\Ext(X, C(K))=0$. The characterization is given
in terms of properties of the metric projection from the space
$Z(X, \R)$ of $\R$-valued $z$-linear maps defined on $X$ onto the
algebraic dual $X'$.

Previous results in this direction were obtained by Kalton and
Pe\l czy\'{n}ski (see \cite{kaltpelc}) who had proved that if
$\Ext(X, C(K))=0$ then $X$ must have the Schur property; then
Kalton shows in \cite{kaltext} that if $\Ext(X,C(K))=0$ then $X$
must also have the strong-Schur property. In \cite{johnzipp},
Johnson and Zippin show that $\Ext(X,C(K))=0$ for every space $X$
dual of a subspace of $c_0$, while Kalton proves in \cite{kaltext}
that the converse is also true when $X$ has an unconditional FDD.
The paper \cite{ccky} gives the first steps towards a theory of
extensions with $C(K)$-spaces by establishing conditions under
which $\Ext(X,C[0,1])\neq 0$ or $\Ext(X,C(\omega^\omega)) \neq
0$.\\

\section{Background on the $z$-linear representation of extensions}
An exact sequence of Banach spaces is a diagram $0\to
Y\stackrel{j}\to E \stackrel{q}\to X\to 0$ composed with Banach
spaces and operators in such a way that the kernel of each arrow
coincides with the image of the preceding. It is sometimes called
an extension of $X$ by $Y$ or, simply, an extension of $X$. The
open mapping theorem makes $Y$ a subspace of $E$ through the
embedding $j$ and $X$ the corresponding quotient space through
$q$.

Following the theory developed by Kalton \cite{kalt} and Kalton
and Peck \cite{kaltpeck} in which extensions $0\to Y\to E \to X\to
0$ of quasi-Banach spaces were identified with quasi-linear maps
$F: X \to Y$, extensions of Banach spaces can be represented (see
\cite{castgonz,cabecastdu}) with a particular type of quasi-linear
maps called $z$-linear maps. These are homogeneous maps such that
for some constant $C$ and every finite set of points $\{z_1,
\dots, z_n \} \subset X$ one has $$\|F (\sum_{i=1}^nz_i
)-\sum_{i=1}^nF(z_i) \| \leq C \sum_{i=1}^n\|z_i\|.$$ The infimum
of the constants $C$ above is denoted $Z(F)$.

We need to introduce some background about $z$-linear maps since
most of our work shall be developed working with this
representation of exact sequences. A $z$-linear map $F: X \to Y$
determines a quasi-norm on the product space $Y \times X$ given by
$\|(y,x)\|_F = \|y - Fx\| + \|x\|$. Let us call this quasi-Banach
space $Y \oplus_F X$. If $co(Y \oplus_F X)$ denotes its Banach
envelope (i.e., the Banach space having as unit ball the closed
convex hull of the unit ball of $\|\cdot \|_F$, it is easy to see
that $co(Y \oplus_F X)$ is $Z(F)$-isomorphic to $Y \oplus_F X$. In
this way each $z$-linear map induces an exact sequence $0 \to Y
\stackrel{j_F}\to co(Y \oplus_F X) \stackrel{q_F}\to X \to 0$ of
Banach spaces with embedding $j_F(y) = (y,0)$ and quotient map
$q_F(y,x) =x$. To obtain a $z$-linear map associated to a given
exact sequence $0\to Y\to E \to X\to 0$ of Banach spaces one can
proceed as follows: take an homogeneous and bounded selection $b:
X\to E $ for the quotient map, then take a linear selection $l: X
\to E$ for the quotient map, and finally make the difference $F =
b-l$, which is a $z$-linear map $X \to Y$. The exact sequences
$0\to Y\to E \to X\to 0$ and $0 \to Y \to co(Y \oplus_F X) \to X
\to 0$ are \emph{equivalent}, in the sense that there exists a
commutative diagram $$
\begin{CD}0@>>> Y@>>> E@>>> X@>>>0\\
&& \Vert&& @VTVV \Vert\\ 0@>>> Y@>>> co(Y \oplus_F X)@>>> X@>>>0.
\end{CD}
$$The classical 3-lemma from homological algebra plus the open
mapping theorem imply that $T$ is an isomorphism. With this
meaning we shall say that the $z$-linear map $F$ has been
associated to the exact sequence and we will represent this using
the notation $0\to Y \stackrel{j}\to E \stackrel{q}\to X\to 0
\equiv F$.\smallskip

Two $z$-linear maps $F$ and $G$ are said to be {\it equivalent},
and written $F\equiv G$, when they induce equivalent exact
sequences. The classical theory establishes that this happens if
and only if there is a linear map $L: X \to Y$ such that $\|F - G
- L\| = \sup \{\|(F -G-L) (x)\| : \|x\|\leq 1 \}< + \infty$. We
will say sometimes that $G$ is a version of $F$. An exact sequence
$0\to Y\stackrel{j}\to E \stackrel{q}\to X\to 0 \equiv F$ is said
to split if it is equivalent to the \emph{trivial} sequence $0\to
Y \to Y\oplus X \to X \to 0$; equivalently, $j(Y)$ is complemented
in $E$. In $z$-linear terms this means $F\equiv 0$, namely, there
is a linear map $L: X \to Y$ such that $\|F - L \|<+\infty$. We
write $\Ext(X,Y)=0$ to mean that every exact sequence with $Y$ as
subspace and $X$ as quotient splits. The lower sequence in a
diagram $$
\begin{CD}0@>>> Y@>j>> X@>>> Z@>>>0\equiv F\\
&& @VTVV@VVV \Vert\\ 0@>>> E@>>> X'@>>> Z@>>>0
\end{CD}$$is called the {\it push-out sequence} and has $TF$ as associated
$z$-linear map. One has that $TF \equiv 0$ if and only if $T$
extends to $X$ through $j$. Dually, the  lower sequence in a
diagram $$
\begin{CD}0@>>> Y@>>> X@>q>> Z@>>>0\equiv F\\
&& \Vert @AAA @AASA\\ 0@>>> Y@>>> X'@>>> E@>>>0
\end{CD}$$is called the {\it pull-back sequence} and has $FS$ as associated
$z$-linear map. One has that $FS \equiv 0$ if and only if $S$ can
be lifted to to $X$ through $q$.

\section{Linearization and factorization of $z$-linear maps}

It will be necessary for us to deal with the space $Z(X,Y)$ of
$z$-linear maps $F: X  \to Y$ considered as mere functions (i.e.,
without equivalence relation). The space $Z(X, \R)$ admits a
semi-normed (not necessarily Hausdorff) topology induced by the
seminorm  $Z(\cdot)$ (the constant of $z$-linearity).

In order to get our main result, Theorem \ref{w*metric}, we need a
refinement of one of the main results in \cite[Thm. 2.1]{ccky}.
There it was shown that for each separable space $X$ there exists
a universal exact sequence $0 \to C[0,1] \to E \to X \to 0 \equiv
U$ with the property that each exact sequence $0 \to C[0,1] \to
E_1 \to X \to 0 \equiv F$ is a push-out of $U$. In particular,
there exists an operator $\phi: C[0,1] \to C[0,1]$ such that $F
\equiv \phi U$. Our aim is to get an equality in the previous
result instead of the mere equivalence. With that purpose in mind
we introduce the linearization and factorization processes, and
recall Zippin's extension method.\\

\noindent \textbf{Linearization process.} Given a $z$-linear map
$F: X \to Y$ and a Hamel basis $(e_\gamma)$ for $Z$, we define a
linear map $\ell_F: X \to Y$ by setting
$\ell_F(e_\gamma)=F(e_\gamma)$. The process $F \to \ell_F$ is
linear. The \emph{linearized form} of $F$ (with respect to a given
Hamel basis) is its version $F-\ell_F$. We shall call the
correspondence $L: Z(X,Y) \to Z(X,Y)$ given by $L(F)  = F -
\ell_F$ \emph{linearization process} (with respect to a given
Hamel basis). We shall omit from now own the coda "with respect to
a given Hamel basis".\\

\noindent \textbf{Factorization process.} Consider a minimal (in
the purely algebraic sense) set $S^+$ such that the unit sphere
$S$ of $X$ coincides with $\cup_{z \in S^+}\{z, -z\}$. Naming
$(e_z)_{z \in S^+}$ the canonical basis of $l_1(S^+)$, we
construct the quotient map $q: l_1(S^+) \to Z$ defined by $q(e_z)
= z$. An homogeneous bounded section for $q$ comes defined as: If
$p \in S^+$ then $b(p) = e_p$; the map $b$ can be extended by
homogeneity to $S$ and then to the whole $X$. The $z$-linear map
$P = b-\ell_b$ is associated with the exact sequence $0 \to K(X)
\to l_1(S^+) \stackrel{q}\to X \to 0$, which we call projective
presentation of $X$. Given a $z$-linear map $F: X \to Y$ we
construct the associated sequence $0 \to Y \stackrel{j_F}\to co(Y
\oplus_F X) \stackrel{q_F}\to X \to 0 \equiv F$. A homogeneous
bounded selection for $q_F$ is $B(x) = (Fx,x)$; it has norm $1$
since $\|Bx\| = \|x\|_F  \leq \|x\|$. A linear selection for $q_F$
is $\Upsilon(x)=(0,x)$ and one has $F= B - \Upsilon$. We define a
norm one operator $\phi(F): l_1(S^+) \to co(Y \oplus_F X)$ as:
$$\phi(F)(\sum_{z \in S^+} \lambda_z e_z ) = \sum_{z \in S^+}
\lambda_z Bqe_z.$$

The restriction of this operator to $K(X)$ shall be called
$\phi_F$, and it takes values in $Y$. Since the norm of $Y$ is
$Z(F)$-equivalent to that of $co(Y \oplus_F X)$ we get $\|\phi_F:
K(X) \to Y \| \leq Z(F)$. We call \emph{factorization process} the
correspondence $Z(X,Y) \to \mathcal L(K(X), Y)$ given by $F \to
\phi_F$ . Observe now that $B = \phi(F) b$ because if $p \in S^+$
then $\phi(F) b(p) = \phi(F)(e_p) = Bqe_p = Bp$. Hence$$\phi_F P =
\phi(F) b - \phi(F) \ell_b = B - \phi(F) \ell_b = B - \Upsilon +
\Upsilon -\phi(F)\ell_b = F + \Upsilon -\phi(F)\ell_b.$$ Since $P$
is linearized form so must be $\phi(F)P$ and therefore $\Upsilon -
\phi(F)\ell_b = -\ell_F$. This means that when the factorization
process is applied to a $z$-linear map in linearized form $G = F -
\ell_F$ then one obtains $\phi_G P =G$.\\

\noindent \textbf{Zippin's extension method.} Let $\delta_X : X
\to C(B_{X^*}, w^*)$ be the canonical embedding. As it was
observed by Zippin in \cite{zippill,zippln}, the $w^*$-continuous
map $\tau: B_{X^*} \to B_{C(B_{X^*})^*}$ defined by
$\tau(x^*)(f)=f(x^*)$ provides an extension for every norm-one
operator $T:X\to C(K)$ to $C(B_{X^*})$ through $\delta_X$ in the
following way: $\widehat{T}(f)(k)=\tau(T^*(\delta_k))(f)$. We will
say that $\widehat{T}$ is the Zippin extension of $T$.\smallskip

Putting together the three processes we get.

\begin{lemma} For every Banach space $X$ there exists a compact space $\Xi[X]$
and a $z$-linear map $\Delta: X \to C(\Xi[X])$ such that for every
$z$-linear map $F: X \to C(K)$ in linearized form there exists an
operator $Z(F)\Phi_F: C(\Xi[X]) \to C(K)$ with norm $\|\Phi_F\|
\leq Z(F)$ such that $Z(F)\Phi_F \Delta = F$.
\end{lemma}
\begin{proof} Consider the extension $0 \to K(X) \to l_1(S^+) \to X \to 0\equiv
P$ and the operator $\phi_F$ given by the factorization process as
described above. Let $\Phi_F$ be Zippin's norm one extension of
$Z(F)^{-1}{\phi_F}$. It is clear that $Z(F)\Phi_F \Delta = Z(F)
\Phi_F \delta_X P = Z(F) Z(F)^{-1}\phi_F P = F$.
\end{proof}

\section{Characterization of the spaces $X$ such that $\Ext(X, C(K))=0$}

Besides the semi-normed topology, we will need to consider on the
space $Z(X, \R)$ the topology $w^*$ of pointwise convergence: we
shall say that $F=w^*-\lim F_\alpha$ if for every $x\in X$ one has
$F(x) = \lim F_\alpha(x)$. A map $Z(X, \R) \to Z(X,\R)$ will be
called $w^*$-continuous if it is $w^*$-continuous on the unit
ball; namely, it transforms $w^*$-convergent nets on the unit ball
into $w^*$-convergent nets. If $X'$ denotes the algebraic dual of
$X$, then $X' \subset Z(X,\R)$ and the restriction of the
$w^*$-topology to $X'$ is the weak $w(X', X)$-topology.\smallskip

The so-called ``nonlinear Hahn-Banach theorem" shown in
\cite{cabecastdu} asserts that given a $z$-linear map $F: X \to
\R$ there exists a linear map $L \in X'$ such that $\|F - L\| \leq
Z(F).$ This makes non-vacuous the following definition:

\noindent \textbf{Definitions.} A map $m: Z(X,\R) \to X'$ will be
called a $\lambda$-metric projection if $$\|F - m(F)\| \leq
\lambda Z(F).$$ If we do not need to emphasize $\lambda$ we just
speak of metric projection. We shall say that a Banach space $X$
admits a metric projection with a given property $P$ if there is a
metric projection $m: Z(X,\R) \to X'$ with property $P$.\\

For instance, every Banach space admits a $1$-metric projection.
In \cite{metric} it was shown that a Banach space $X$ is an
$\mathcal L_1$-space if and only if it admits a {\it linear}
metric projection, which  This yields that $X$ admits  a {\it
linear} metric projection if and only if $\Ext(X, Y^*)=0$ for
every dual space.  A characterization of Banach spaces such that
all their extensions by any $C(K)$-space are trivial can also be
obtained in terms of properties of metric projections.

\begin{thm}\label{w*metric} A Banach space $X$ admits a
$w^*$-continuous metric projection if and only if  $\Ext(X,
C(K))=0$ for every compact Hausdorff space $K$.
\end{thm}
\noindent{\bf Proof of the necessity.} Let $F: X \to C(K)$ be a
$z$-linear map and let $m:  Z(X,\R) \to X'$ be a $w^*$-continuous
$\lambda$-metric projection. We define a linear map $M: X \to
C(K)$ by $$M(x)(k) = m(\delta_k F)(x),$$where $\delta_k$ is the
evaluation at $k$. The map $M$ is well defined since $M(x)$ is a
continuous function: whenever $k = \lim k_\alpha$ on $K$ then
$M(x)(k) = m(\delta_k F)(x) = m(w^*-\lim \delta_{k_\alpha}F) (x) =
\lim m\left( \delta_{k_\alpha}F\right )(x) = \lim M(x)(k_\alpha)$.
Moreover, $|F(x)(k) - M(x)(k)| = |\delta_kF(x) - m(\delta_kF)(x)|
\leq \lambda Z(\delta_kF)\|x\|$, which implies $\|F - M\| \leq
\lambda Z(F)$.\\

\noindent {\bf Proof of the sufficiency .} If $\Ext(X, C(K))=0$
for every $C(K)$-space then $\Ext(X, C(\Xi[X]))=0$. In particular,
$\Delta \equiv 0$ (the "initial`` $z$-linear map of Lemma 1) and
there exists a linear map $\Lambda: X \to C(\Xi[X])$ such that
$\|\Delta - \Lambda\| < +\infty$. Let $F: X \to \R$ be a
$z$-linear map with $Z(F)\leq 1$. Consider the projective
presentation  $0\to K(X)\to l_1(S^+)\to X\to 0\equiv P$ and the
operator $\phi_F: K(X)\to \R$ such that $\phi_FP = F - \ell_F$
given by the factorization process. Let $\Phi_F: C(\Xi[X])\to \R$
be Zippin's extension. We define a map $m: Z(X, \R) \to X'$ as
$$m(F)= \Phi_F \Lambda.$$ To prove that $m(\cdot)$ is
$w^*$-continuous we decompose it in three applications:
\begin{enumerate}
\item The linearization process $L(F)= F - \ell_F$, which is $w^*$-continuous;
\item The factorization process. We show now it is
w*-continuous. Observe that if we restrict ourselves to work with
the subspace $\varphi_0(S^+)$ of $l_1(S^+)$ of all finitely
supported sequences then $ K_0 = K(X) \cap \varphi_0(S^+)$ is
dense in $K(X)$. On this dense subspace the operator $\phi_F$
takes the form $$ \phi_F (\sum\lambda_z e_z) = \sum \lambda_z
Fqe_z.$$ Tus, if $\{G_\alpha\}$ is a net in the unit ball such
that $G = w^*-\lim G_\alpha$ then $\phi_G (u)= \lim
\phi_{G_\alpha(u)}$ for all $u \in K_0$. Since
$\|\phi_{G_\alpha}\| \leq 1$, the sequence of uniformly bounded
operators is $w^*$-convergent on the whole $K(X)$.
\item Zippin's extension method is a $w^*$-continuous process $B_{K(X)^*} \to B_{C(\Xi[X])^*}$
as we have already remarked.
\end{enumerate}

The $w^*$-continuous metric projection we are looking for is$$F
\to m(F) + \ell_F.$$ It is obviously $w^*$-continuous and it
remains to show that it is a metric projection:$$ \| F - m(F) -
\ell_F \| = \| \Phi_F \Delta - \Phi_F\Lambda \| \leq \| \phi_F\|\|
\Delta - \Lambda\| \leq Z(F) \| \Delta - \Lambda\|.$$

\noindent \textbf{Remark.} In quantitative terms the previous
arguments provide that every $C(K)$-valued $z$-linear map admits a
linear map at distance at most $\lambda Z(\cdot)$ if and only if
there is a $w^*$-continuous method to assign to each $\R$-valued
$z$-linear map a linear map at distance at most $\lambda
Z(\cdot)$.

\section{Remarks and open questions}

1. \emph{Existence of w*-continuous metric projections.} In
general, every Banach space $X$ admits both $1$-metric projections
(such as that given by the non-linear Hahn-Banach theorem) and
$w^*$-continuous projections (such as the linear map $F \to
\ell_F$). As we have already seen, it is extremely unusual that
the former are $w^*$-continuous and the latter is (usually) not a
metric projection.\smallskip

2. \emph{Existence of w*-continuous $C$-metric projections on
finite dimensional spaces.} There is one instance in which the the
linear map $F \to \ell_F$ is a metric projection: when the space
is $l_1$, we work on the dense subspace of finitely supported
sequences and take as Hamel basis the canonical Schauder basis of
$l_1$. The final step of passing from a dense subspace to the
whole space is a classical extension result for $z$-linear maps
(see \cite{kaltpeck,castgonz}). It is therefore clear that
finite-dimensional spaces $E$ admit a $w^*$-continuous
$\mathrm{dist} (E, l_1^{dim E})$-metric projection. However, as
the dimension of $E$ increases the "$\lambda$-metric" character of
the projection is spoiled. Nevertheless, it is possible to get
$w^*$-continuous $\lambda$-metric projections on finite
dimensional spaces with uniform constant $\lambda$ independent of
the dimension. To do that observe that $\mathcal L_{\infty,
\mu}$-spaces are \emph{locally complemented} in any superspace,
which means (see \cite{kaltloc}) that there exists a constant
$C(\mu)$ such that for every finite-dimensional space $E$ every
$z$-linear map $F: E \to \mathcal L_{\infty, \mu}$ admits a linear
map $\ell: E \to \mathcal L_{\infty, \mu}$ such that $\|F - \ell\|
\leq C(\mu)$. Since $C(K)$-spaces are $\mathcal L_{\infty, 1+
\varepsilon}$-spaces, it turns out that finite dimensional spaces
admit $w^*$-continuous $C(1+\varepsilon)$-metric
projections.\smallskip

3. \emph{Existence of $w^*$-continuous metric projections in
$l_1(X_n)$.} To work with infinite dimensional spaces we recall
again the extension result for $z$-linear maps; thus, we only need
to work on a dense subspace; i.e., if $X_0$ denotes a dense
subspace of $X$, then $X_0$ admits a $w^*$-continuous metric
projection if and only if $X$ does.

\begin{proposition}\label{amalgama} If $A_n$ is a sequence of spaces
admitting $w^*$-continuous $\lambda$-metric projections then the
vector sum $l_1(A_n)$ admits a $w^*$-continuous
$(1+\lambda)$-metric projection.
\end{proposition}
\begin{proof} Let $\varphi_0(A_n)$ be the dense subspace of $l_1(A_n)$
formed by all finitely supported sequences. The map
$Z(\varphi_0(A_n), \R) \to l_\infty( Z(A_n,\R))$ that assigns to a
$z$-linear map $F$ the sequence $(F_{|A_n})_n$ of its restrictions
is clearly $w^*$-continuous. If $m_n :  Z(A_n, \R) \to A_n'$ are
$w^*$-continuous metric projections then we can define a
$w^*$-continuous metric projection $m: Z(\varphi_0(A_n), \R) \to
l_1(A_n)' $ as follows: take $(e_\gamma^n)_\gamma$ a Hamel basis
for $A_n$, and set $m(F)(e_\gamma^n) = m_n(F_n)(e_\gamma^n).$ If
$a=(a_n) \in \varphi_0(A_n)$ then $$|Fa - m(F)a| = |Fa - \sum
F_na_n + \sum F_na_n - \sum m_n(F_n)a_n| \leq (1+\lambda)Z(F)\sum
\|a_n\|.$$
\end{proof}

4. \emph{The role of the $C(K)$-space.} It is not difficult to see
that for separable spaces $X$ the condition ``$\Ext(X, C(K))=0$
for all compact spaces $K$" is equivalent to "$\Ext(X,
C[0,1])=0$``. It has been shown in \cite{ccky} that $\Ext(X,
C(\omega^\omega))=0$ is a strictly weaker condition; precisely, if
$T$ denotes the dual of the original Tsirelson space then $\Ext(T,
C[0,1])\neq 0$ while $\Ext(T, C(\omega^\omega))=0$. It seems a
touchy question how metric projections would reflect the change of
the $C(K)$ target space. For instance, it is clear that, without
separability assumptions, if $X$ admits a metric projection which
is $w^*$-continuous at $0$ then $\Ext(X, c_0)=0$.

\textbf{Question.} Is it true that $\Ext(X, c_0)=0$ if and only if
$X$ admits a metric projection $w^*$-continuous at $0$? In
particular: Does every separable Banach space admit a metric
projection $w^*$-continuous at $0$?\smallskip

6. \emph{On the equation $\Ext(X, C(K))=0$ for subspaces of
$l_1(\Gamma)$.} Returning to Kalton's proposition that separable
Banach spaces such that $\Ext(X, C[0,1])=0$ must have the
strong-Schur property, it makes sense to approach the converse by
asking:

\textbf{Question.} Is $\Ext(M, C[0,1])=0$ for every subspace $M$
of $l_1$?\smallskip

Let us show the existence of a subspace $K(Z)$ of $l_1(\Gamma)$
such that $\Ext(K(Z), c_0) \neq 0$, which implies that $K(Z)$ does
not admit a metric projection $w^*$-continuous at $0$. To get the
example, recall that $l_\infty/c_0$ is not injective (see
\cite{rosestu}) and thus there exists some (necessarily
nonseparable) space $Z$ for which $\Ext(Z, l_\infty/c_0)\neq 0$.
This gives a nontrivial exact sequence $0 \to l_\infty/c_0 \to E
\to Z \to 0 \equiv F$. Given an exact sequence $0 \to K(Z)
\stackrel{i}\to l_1(\Gamma) \to Z \to 0 \equiv P$ there is an
operator $\phi: K(Z) \to l_\infty/c_0$ such that $F \equiv \phi
P$. Consider the exact sequence $0 \to c_0 \to l_\infty
\stackrel{q}\to l_\infty/c_0 \to 0 \equiv I$. If $I\phi \equiv 0$
then $\phi$ can be lifted to an operator $\eta: K(Z) \to l_\infty$
through $q$; this operator $\eta$, in turn can be extended to an
operator $\psi: l_1(\Gamma) \to l_\infty$ through $i$. The
operator $q\psi: l_1(\Gamma) \to l_\infty/c_0$ extends $\phi$
since $q\psi i= q\eta = \phi$. Therefore $F \equiv \phi P \equiv
0$, against the hypothesis.\smallskip

We close the paper observing that a remarkable class of spaces
with the strong-Schur property are those having a $l_1$-skipped
blocking decomposition introduced by Bourgain and Rosenthal
\cite{bourrose}. It would be interesting to know if those spaces
admit $w^*$-continuous metric projections.

\end{document}